\documentclass[11pt,a4paper]{article}
\usepackage[utf8]{inputenc}
\usepackage{epic,eepic,epsfig,amsmath,amssymb,amsthm}
\usepackage{t1enc}

\usepackage[francais,english]{babel}

\usepackage{caption}

\usepackage{color}

\usepackage{bbm}

\def\virgp{\raise 2pt\hbox{,}}

\renewcommand{\geq}{\geqslant}
\renewcommand{\leq}{\leqslant}
\def\N{{\mathbb N}}

\def\R{{\mathbb R}}

\def\virgp{\raise 2pt\hbox{,}}
\def\cdotpv{\raise 2pt\hbox{;}}

\def\1{\mathbbm{1}}

\newtheorem{theorem}{Theorem}[section]

\newtheorem{proposition}[theorem]{Proposition}

\theoremstyle{remark}
\newtheorem{remark}{Remark}[section]

\theoremstyle{definition}
\newtheorem{definition}{Definition}[section]

\theoremstyle{definition}

\theoremstyle{definition}

\begin{document}
	
	\title{The Finite Difference Method,\\ for the heat equation on Sierpi\'{n}ski simplices}
	
	\author{Nizare Riane, Claire David\footnote{Corresponding author: Claire.David@Sorbonne-Universite.fr}}

	\maketitle
	\centerline{Sorbonne Universit\'e}
	
	\centerline{CNRS, UMR 7598, Laboratoire Jacques-Louis Lions, 4, place Jussieu 75005, Paris, France}

	\begin{abstract}
		
		In the sequel, we extend our previous work on the Minkowski Curve to Sierpi\'{n}ski simplices (Gasket and Tetrahedron), in the case of the heat equation. First, we build the finite difference scheme. Then, we give a theoretical study of the error, compute the scheme error, give stability conditions, and prove the convergence of the scheme. Contrary to existing work, we do not call for approximations of the eigenvalues.
	\end{abstract}

	\maketitle
	\vskip 1cm
	

	\vskip 1cm
	
	\vskip 1cm
	
	\vskip 1cm
	
	\section{Introduction}

	\hskip 0.5cm Following the seminal work of J.~Kigami~\cite{Kigami1989},~\cite{Kigami1993},
	~\cite{Kigami2001},~\cite{Kigami2003} in the field of analysis on fractals, the natural step was to explore the numerical aspects. \\
	
	It has been initiated, in the case of the Sierpi\'{n}ski gasket~$\mathcal{SG}$, by~K.~Dalrymple, R.~S.~Strichartz, and J.~Vinson~\cite{DalrympleFDM}, who gave an equivalent method for the finite difference approximation. More precisely, the authors use the spectral shape of the solution (heat kernel), which involves eigenvalues and eigenvectors, and therefore calls for an approximation of the eigenvalues. This work has been followed by the one of M.~Gibbons, A.~Raj and R.~S.~Strichartz~\cite{GibbonsFEM}, where they describe how one can build approximate solutions, by means of piecewise harmonic, or biharmonic, splines, again in the case of~$\\mathcal{SG}$.  They go so far as giving theoretical error estimates, through a comparison with experimental numerical data. \\

	After our work~\cite{RianeDavidM}, where we built a Laplacian on the Minkowski curve, we went so far as implementing the resulting finite difference scheme, for which we could not find any equivalent in the existing literature.  \\

	The novelty of our contribution lies in defining the discretization of the considered PDEs (heat and wave equations), by taking into account the recursive construction of the matrix related to the sequence of graph Laplacians. Contrary to the aforementioned work, we thus did not call for approximations of the eigenvalues. This enabled us not only to compute the consistency error, but, also, to set stability conditions of Courant-Friedrichs-Lewy type, and, then, to prove the  convergence of the scheme.\\
	
	In the sequel, we extend this method to Sierpi\'{n}ski simplices (Gasket and Tetrahedron), in the case of the heat equation. First, we build the finite difference scheme. Then, we give a theoretical study of the error, compute the scheme error, give stability conditions, and prove the convergence of the scheme, both for the explicit and the implicit Euler methods. For the explicit scheme, we obtain a necessary and sufficient CFL condition of the form~$h\lesssim(\delta x)^{D_W}$, where~$\delta x$ denotes the spatial step and~$D_W$ the walk dimension of the simplex, which replaces the exponent~$2$ of the smooth case. The implicit scheme is unconditionally stable. In both cases, the error is of order~$\mathcal{O}(h)+\mathcal{O}\left((\delta x)^{\frac{D_W-D_H}{2}}\right)$, where~$D_H$ denotes the self-similar dimension.
	
	\vskip 1cm

	\section{The Sierpi\'{n}ski simplices}
	
	In the sequel, we place ourselves in the Euclidean space of dimension~$d-1$ for a strictly positive integer~$d$, referred to a direct orthonormal frame. The usual Cartesian coordinates will be denoted by~$(x_1,x_2,...,x_{d-1})$.\\
	
	\noindent Let us introduce the family of contractions~$f_i$,~\mbox{$1\leq i \leq d$}, of fixed point~$P_{i-1}$ such that, for any~$X\,\in\,\R^{d-1}$, and any integer~$i$ belonging to~\mbox{$\left \lbrace 1, \hdots,d \right \rbrace $}:
	
	$$f_i(X)=\frac{1}{2}(X+P_{i-1})$$

	\vskip 1cm
	
	 	\noindent According to~\cite{Hutchinson1981}, there exists a unique nonempty compact subset $\mathfrak{SS} \subset \R^{d-1}$ such that:
 
	\[\mathfrak{SS} = \underset{  i=1}{\overset{d}{\bigcup}}\, f_i(\mathfrak{SS})\]
	\noindent which will be called the Sierpi\'{n}ski simplex.\\

	\vskip 1cm

	\noindent We will denote by~$V_0$ the ordered set of the points:
	
	$$\left \lbrace P_{0},\hdots,P_{d-1}\right \rbrace$$

	 	\noindent The set of points~$V_0$, where, for any~$i$ of~\mbox{$\left \lbrace  0,...,d-1  \right \rbrace$}, every point~$P_i$ is linked to the others, constitutes a complete graph, that we will denote by~$ {\mathfrak {SS}}_0$.~$V_0$ is called the set of vertices of the graph~$ {\mathfrak {SS}}_0$.\\

	\noindent For any strictly positive integer~$m$, we set:
	$$V_m =\underset{  i=1}{\overset{d}{\bigcup}}\, f_i \left (V_{m-1}\right )$$

	\noindent The set of points~$V_m$, where the points of an~$m^{th}$-order cell are linked in the same way as ${\mathfrak {SS}}_0$, is a graph, which we will denote by~$ {\mathfrak {SS}}_m$.~$V_m$ is called the set of vertices of the graph~$ {\mathfrak {SS}}_m$. We will denote, in the following, by~${\mathcal N}_m$ the number of vertices of the graph~$ {\mathfrak {SS}}_m$.

	\vskip 1cm

	\begin{proposition}
		Given a natural integer~$m$, we will denote by~$\mathcal{N}_m$ the number of vertices of the graph ${\mathfrak {SS}}_m$. One has:
		
		$$\mathcal{N}_0  =d$$
		
		\noindent and, for any strictly positive integer~$m$:
		
		$$\mathcal{N}_m   =d\,\mathcal{N}_{m-1}- \frac{d\,(d-1)}{2}$$
	\end{proposition}

	\vskip 1cm
	
	\begin{proof} 	The graph $\mathfrak{SS}_m$ is the union of~$d$ copies of the graph $\mathfrak{SS}_{m-1}$. Each copy shares exactly one vertex with each one of the other copies. So, one may consider the copies as the vertices of a complete graph~$K_d$, the number of edges of which is equal to~$\displaystyle\frac{d\,(d-1)}{2}$, which leads to~$\displaystyle\frac{d\,(d-1)}{2}$ vertices counted twice.
		
	\end{proof}
	
	\vskip 1cm
	
	\begin{remark}
		\noindent One may check that~$\mathcal{N}_m=\displaystyle\frac{d^{m+1}+d}{2}$.
	\end{remark}
	
	\vskip 1cm
	
	\section{The finite difference method on the Sierpi\'{n}ski simplices}

	In the sequel, we will denote by~$T$ a strictly positive real number, by~$\mathcal{N}_0$ the cardinal of $V_0$, and by~$\mathcal{N}_m$ the cardinal of $V_m$.\\

	\subsection{The heat equation}
	
	\subsubsection{Formulation of the problem}

	\noindent We may now consider a solution~$u$ of the problem:
	
	$$ \left \lbrace  \begin{array}{ccccc}
	\displaystyle \frac{\partial u}{\partial t}(t,X)-\Delta u(t,X)&=&0 & \forall \, (t,X)\,\in \,\left]0,T\right[  \times {\mathfrak{SS}}\\
	u(t,X)&=&0 & \forall \, ( t,X) \,\in \,\left[0,T\right[ \times \partial {\mathfrak{SS} }\\
	u(0,X)&=&g(X)  & \forall \,X\,\in \,{\mathfrak{SS}}
	\end{array}\right.$$
	
	\noindent where~$\Delta$ denotes the Laplacian on~$\mathfrak{SS}$ with respect to the standard self-similar measure~$\mu$, which gives the weight~$d^{-m}$ to each~$m^{th}$-order cell, and where~$\partial\mathfrak{SS}=V_0$.\\
	
	\noindent In order to define a numerical scheme, one may use a first order forward difference scheme to approximate the time derivative~$\displaystyle \frac{\partial u}{\partial t}$. The Laplacian is approximated by means of the graph Laplacians~$\Delta_m \,u$, defined on the sequence of graphs~$\left ( {\mathfrak{SS}}_m\right)_{m\in\N^{\star}}$. \\
	
	\noindent To this purpose, we fix a strictly positive integer~$N $, and set:
	
	$$ \displaystyle{h =\frac{T}{N}}$$
	
	\noindent One has, for any integer~$k$ belonging to~$\left \lbrace 0, \hdots, N-1 \right \rbrace$, and for~$u$ twice continuously differentiable with respect to the time variable:
	
	$$   \forall \,X\,\in \, {\mathfrak{SS}} \,:   \quad
	\displaystyle\frac{\partial u}{\partial t}(kh,X)=\displaystyle\frac{1}{h}\,\left( u((k+1)\, h,X)-u(kh,X)\right)+{\mathcal O}(h)
	$$
	
	\noindent According to~\cite{StrichartzLivre2006}, the Laplacian on Sierpi\'{n}ski simplices~$\mathfrak{SS}$ is given by:
	
	$$ \forall \, X\,\in \,{\mathfrak{SS}}\setminus V_0 \,:   \quad
	\Delta u(t,X)= 	\displaystyle\lim_{m\rightarrow +\infty} r^{-m}\left(\int_{\mathfrak{SS}} \psi^{(m)}_{X_m} \, d\mu\right)^{-1} \, \sum_{Y \underset{m}\sim X_m} \left( u(t,Y)- u(t,X_m)\right)
	$$
	
	\noindent where, for any natural integer~$m$,~\mbox{$X_m\, \in\, V_m\setminus V_0$},~$\psi^{(m)}_{X_m}$ denotes the piecewise harmonic function (harmonic spline) taking the value~$1$ on~$X_m$ and~$0$ on the other points of~$V_m$,~$r$ denotes the energy renormalization ratio, and:

	$$	\displaystyle\lim_{m\rightarrow \infty}X_m=X$$
	
	\noindent In the case of Sierpi\'{n}ski simplices, one has (see~\cite{Fukushima1992},~\cite{StrichartzLivre2006}):
	
	$$r=\frac{d}{d+2} \quad , \quad \int_{\mathfrak{SS}} \psi^{(m)}_{X_m} \, d\mu=\frac{2}{d^{m+1}}$$
	
	\noindent since the support of~$\psi^{(m)}_{X_m}$ consists of the two~$m^{th}$-order cells containing~$X_m$. This yields:
	
	$$r^{-m}\left(\int_{\mathfrak{SS}} \psi^{(m)}_{X_m} \, d\mu\right)^{-1}=\frac{d}{2}\,(d+2)^m$$
	
	\noindent This enables one to approximate the Laplacian, at a~$m^{th}$ order,~$m \,\in\,\N^\star$, using the graph normalized Laplacian as follows:
	
	$$ \forall \, k\,\in\, \left \lbrace 0, \hdots, N-1 \right \rbrace,\,\forall \,X\,\in \,V_m\setminus V_0  \,: \quad
	\Delta u(kh,X)\approx \frac{d}{2}\,(d+2)^m \,\displaystyle\sum_{Y \underset{m}\sim X} \left( u(kh,Y)- u(k\, h,X)\right)
	$$
	
	\noindent By combining those two relations, one gets the following scheme, for any integer~$k$ belonging to~$\left \lbrace 0, \hdots, N-1 \right \rbrace$, any point~$P_j$ of~$V_0$,~$0\leq j \leq d-1$, and any~$X$ in the set~\mbox{$V_m\setminus V_0$} :
	
	$$\left ({\mathcal S}_{\mathcal H}\right) \quad \left \lbrace \begin{array}{cccc}
	\displaystyle \frac{u^m_h((k+1)\,h,X)-u^m_h(k\,h,X)}{h}&=& \displaystyle\frac{d}{2}\,(d+2)^m \,\displaystyle \sum_{Y \underset{m}\sim X} \left( u^m_h(k\,h,Y)- u^m_h(kh,X)\right) &   \\
	u^m_h(k\,h,P_j)&=&0 &     \\
	u^m_h(0,X)&=&g(X) &
	\end{array} \right.$$
	
	\noindent Let us define the approximate equation as:
	
	$$
	u^m_h((k+1)\,h,X)=u^m_h(k\, h,X) + h\, \frac{d}{2}\,(d+2)^m \, \sum_{Y \underset{m}\sim X} \left( u^m_h(k\, h,Y)- u^m_h(k\,h,X)\right) \quad
	\forall \, k\, \in\,\left \lbrace 0, \hdots, N-1 \right \rbrace,  \, \forall \,X\,\in \,V_m\setminus V_0
	$$
	
	\noindent We now fix~$m\,\in \,\N$, and denote any~$X\,\in \,V_m\setminus V_0$ as~$X_{w,P_i}=f_w(P_i)$, where~\mbox{$w\,\in\,\{1,\dots,d\}^m$} is a word of length $m$,~$f_w=f_{w_1}\circ\hdots\circ f_{w_m}$, and where $P_i$,~$0 \leq i \leq d-1$, belongs to~$V_0$. Since a point may have several addresses, we assign to each point its address of lowest lexicographic order. This induces a natural ordering~$X_1,\hdots,X_{\mathcal{N}_m-d}$ of the points of~$V_m\setminus V_0$.\\
	
	\noindent This enables one to introduce, for any integer~$k$ belonging to~$ \left \lbrace 0, \hdots, N \right \rbrace$, the solution vector $U(k)$ as:
	
	\begin{align*}
	U^m_h(k) &=\left(
	\begin{matrix}
	u^m_h(k\,h,X_1)\\
	\vdots\\
	u^m_h(k\,h,X_{\mathcal{N}_m-d})\\
	\end{matrix}
	\right)\\
	\end{align*}
	
	\noindent which satisfies the recurrence relation:
	
	$$U^m_h(k+1) = A \, U^m_h(k)$$
	
	\noindent where:
	
	$$
	A=I_{\mathcal{N}_m -d}-h\,\tilde{\Delta}_{m}
	$$
	
	\noindent  and where~$I_{\mathcal{N}_m -d}$ denotes the~$ {(\mathcal{N}_m -d)}\times  {(\mathcal{N}_m -d)}$ identity matrix, and $\tilde{\Delta}_{m}$ the~$ {(\mathcal{N}_m -d)}\times  {(\mathcal{N}_m -d)}$ normalized Laplacian matrix:
	
	$$\tilde{\Delta}_{m}=\frac{d}{2}\,(d+2)^m\,L_m$$
	
	\noindent where~$L_m$ denotes the matrix of~$-\Delta_m$ restricted to~$V_m\setminus V_0$ (Dirichlet condition): its diagonal coefficients are equal to~$2\,(d-1)$, and, for~$i\neq j$, its coefficient~$(i,j)$ is equal to~$-1$ if~$X_i\underset{m}{\sim}X_j$, and to~$0$ otherwise.
	
	\vskip 1cm
	
	\subsubsection{Theoretical Study of the Error, for H\"{o}lder Continuous Functions}
	
\hskip 0.5cm  In the spirit of the work of~Robert~S.~Strichartz~\cite{Strichartz1999},~\cite{Strichartz2000}, it is interesting to consider the case of H\"{o}lder continuous functions. Why? First, H\"{o}lder continuity implies continuity, which is a required condition for functions in the domain of the Laplacian (we refer to our work~\cite{RianeDavidM} for further details).\\

 Second, a H\"{o}lder condition for such a function will result in fruitful estimates for its Laplacian, which is a limit of difference quotients.\\

	\newpage

  Let us thus consider a function~$u$ in the domain of the Laplacian, and the positive real constant
	
$$\alpha=\displaystyle\frac{\ln(r)}{2\ln(R)}=\frac{\ln\left(\frac{d+2}{d}\right)}{2\,\ln 2}=\frac{D_W-D_H}{2}\,\in\,\left]0,1\right[$$

\noindent where~$r=\displaystyle\frac{d}{d+2}$ is the energy renormalization ratio and~$R=\displaystyle\frac{1}{2}$ the contraction ratio of the similarities~$f_i$. Any function of finite energy on~$\mathfrak{SS}$ is H\"{o}lder continuous of exponent~$\alpha$ with respect to the Euclidean distance (see~\cite{StrichartzLivre2006}). We assume that~$u$ satisfies
	
$$\forall\, (X,Y)\,\in \mathfrak{SS}^2,\forall\,t>0\, : \quad |u(t,X)-u(t,Y)| \leq C(t) \, |X-Y|^\alpha$$
	
\noindent where~$C$ denotes a positive function of the time variable~$t$.\\
	
 Given a strictly positive integer~$m$, we recall that
	
$$\Delta_m u(t,X) = \displaystyle\sum_{Y \in V_m,\,Y\underset{m}{\sim} X} \left (u(t,Y)-u(t,X)\right)  \quad \forall\,t>0,\, \forall\, X\,\in\, V_m\setminus V_0 $$

\noindent and that the Laplacian is defined as the limit 
	
$$\Delta_{\mu} u(t,X)=\displaystyle\lim_{m\rightarrow+\infty} r^{-m} \left( \int_{\mathfrak{SS}} \psi_{X_m}^{(m)} d\mu\right)^{-1}\, \Delta_m u(t, X_m)  \quad \forall\,t>0,\, \forall\, X\,\in\, \mathfrak{SS}\setminus V_0$$
	
\noindent where $\left (X_m\,\in\,V_m\setminus V_0 \right)_{m\in\N}$ is a sequence of points such that:
	$$\displaystyle\lim_{m\rightarrow+\infty}X_m = X$$
	
\noindent and where~$r$ denotes the normalization constant, and~$\psi_{X_m}^{(m)}$ a harmonic spline function.\\
	
  Let us now introduce a strictly positive number~$\delta_{ij}=|P_i-P_j|$, for any~$P_i$ belonging to the set~\mbox{$V_0$}, and any~$P_j$ such that $P_j{\sim} P_i$. On the one hand, we set: $\delta=\max_{i,j}\delta_{ij}$.\\
	
 On the other hand, we define $R$ as the maximum contraction ratio of the similarities~$f_i$ (note that, in the case of Sierpi\'{n}ski simplices, \mbox{$R=\displaystyle\frac{1}{2}$}).\\
	
  For $X\underset{m}{\sim} Y$, $X$ and $Y$ have addresses such that:
	
$$X=f_{\mathcal{W}}(P_i) \quad , \quad Y=f_{\mathcal{W}}(P_j)$$
	
\noindent for some $P_i$ and $P_j$ in $V_0$ and $\mathcal{W}\in \Sigma_m$. May one set 
	
$$R(\mathcal{W}) =R_{\mathcal{W}_1}R_{\mathcal{W}_2} \hdots R_{\mathcal{W}_m}=R^m$$
	
\noindent we then obtain that
	
$$\left |X-Y \right | =\left |f_{\mathcal{W}}(P_i)-f_{\mathcal{W}}(P_j) \right |= R(\mathcal{W}) \left |P_i-P_j \right |  \leq R^m \,\delta   \, \cdot $$
	
  Moreover, any~$X\,\in\,V_m\setminus V_0$ is the junction point of exactly two~$m^{th}$-order cells, and thus has exactly~$2\,(d-1)$ neighbours in~$V_m$.\\
	
  One has then, for any~$X$ belonging to the set~\mbox{$ V_m\setminus V_0$}, any integer~$k$ belonging to~$\left \lbrace 0, \hdots, N-1 \right \rbrace$, and any strictly positive number~$h$,
	
$$\begin{array}{ccc} \left |r^{-m}\, \left( \displaystyle\int_{\mathfrak{SS}} \psi_{X_m}^{(m)} d\mu\right)^{-1} \right | \, h \, \left |\Delta_m u(k\,h,X)\right | &\leq & \left | r^{-m} \, \left( \displaystyle\int_{\mathfrak{SS}} \psi_{X_m}^{(m)} d\mu\right)^{-1} \right |\, h \, \displaystyle\sum_{Y \in V_m,\,Y\underset{m}{\sim} X} \left |u(k\,h,Y)-u(k\,h,X)\right| \\
&\leq & \left | r^{-m} \,\left( \displaystyle\int_{\mathfrak{SS}} \psi_{X_m}^{(m)} d\mu\right)^{-1} \right | \,h\,C(k\,h) \,\displaystyle\sum_{Y \in V_m,\,Y\underset{m}{\sim} X}   |X-Y|^\alpha  \\
&\leq & \left | r^{-m}\,  \left( \displaystyle\int_{\mathfrak{SS}} \psi_{X_m}^{(m)} d\mu\right)^{-1} \right | \,h\,C(k\,h) \,\displaystyle\sum_{Y \in V_m,\,Y\underset{m}{\sim} X}   \displaystyle \delta^{\alpha} \,R^{m\,\alpha}  \\
&=  & 2\,(d-1)\,\delta^{\alpha} \, R^{m\,\alpha}\, \left |r^{-m}\, \left( \displaystyle\int_{\mathfrak{SS}} \psi_{X_m}^{(m)} d\mu\right)^{-1} \right | \,h\,C(k\,h)  \, \cdot \\
\end{array}$$
	
  The increment over one time step of the scheme~\mbox{$\left ({\mathcal S}_{\mathcal H}\right) $}, applied to the function~$u$, thus satisfies
	
$$  \left |r^{-m}\, \left( \displaystyle\int_{\mathfrak{SS}} \psi_{X_m}^{(m)} d\mu\right)^{-1} \right | \, h \, \left |\Delta_m u(k\,h,X)\right | \leq
2\,(d-1)\,\delta^{\alpha} \,  \left(\frac{1}{2}\right)^{m\,\alpha} \left |r^{-m}\, \left( \displaystyle\int_{\mathfrak{SS}} \psi_{X_m}^{(m)} d\mu\right)^{-1} \right | \,h\,C(k\,h) \, \cdot $$
	
\vskip 1cm

  One may note that a sufficient condition for these increments to go to zero is
	
$$ \displaystyle \lim_{m\to + \infty,\, h \to 0^+} \left(\frac{1}{2}\right)^{m\,\alpha}\left | r^{-m} \,\left( \displaystyle\int_{\mathfrak{SS}} \psi_{X_m}^{(m)}\, d\mu\right)^{-1} \right | \,h\,C(k\,h)=0   \, \cdot $$
	
  We take~$C=\sup_k C(k h)$, it reduces to the sufficient condition

\begin{align}
\displaystyle \lim_{m\to + \infty,\, h \to 0^+} \left(\frac{1}{2}\right)^{m\,\alpha}\left | r^{-m} \,\left( \displaystyle\int_{\mathfrak{SS}} \psi_{X_m}^{(m)}\, d\mu\right)^{-1} \right | \,h &=0   \, \cdot
\end{align}

According to \cite{Fukushima1992}, in the case of Sierpi\'{n}ski simplices, $r^m=\displaystyle\left(\frac{d}{d+2}\right)^m$,~$\alpha=\displaystyle-\frac{\ln\left(\frac{d}{d+2}\right)}{2\ln(2)}$, so that~$\displaystyle\left(\frac{1}{2}\right)^{m\,\alpha}=\left(\frac{d}{d+2}\right)^{\frac{m}{2}}$, and a simple calculus as in \cite{StrichartzLivre2006} shows that \mbox{$\displaystyle\int_{\mathfrak{SS}} \psi_{X_m}^{(m)}\, d\mu=\frac{2}{d^{m+1}}$}. The relation reduces to

\begin{align}
\label{eqn:NecessaryFDM}
\displaystyle \lim_{m\to + \infty,\, h \to 0^+} \frac{d\,\left(d\,(d+2)\right)^{\frac{m}{2}}}{2} \,h &=0   \, \cdot
\end{align}

  Since~$\sqrt{d\,(d+2)}<d+2$, this condition is satisfied as soon as~$h\,(d+2)^m$ remains bounded, which is the case under the CFL stability condition~\ref{eqn:ExplicitSpectrum} obtained in the sequel.
	
\vskip 1cm

\subsubsection{Consistency, stability and convergence}

	\paragraph{The scheme error}$\,$\\

	\noindent Let us consider the solution~$u$ of the problem, and assume that~$\displaystyle\frac{\partial^2 u}{\partial t^2}$ is bounded on~$[0,T]\times\mathfrak{SS}$. For all~$k$ in~$\left \lbrace 0, \hdots, N-1 \right \rbrace$ :

	$$   \forall \,X\,\in \, {\mathfrak{SS}} \,:   \quad
	\displaystyle\frac{\partial u}{\partial t}(k\,h,X)=\displaystyle\frac{1}{h}\,\left( u((k+1)\, h,X)-u(k\, h,X)\right)+{\mathcal O}(h)
	$$

	\noindent On the other hand, given a strictly positive integer~$m$,~$X\in V_m\setminus V_0$, and the harmonic spline~$\psi_X^{(m)}$, harmonic on each~$m^{th}$-order cell, taking the value $1$ on $X$ and $0$ on the other vertices of~$V_m$ (see \cite{Strichartz1999}):

	\begin{align*}
	\displaystyle\int_{\mathfrak{SS}} \psi_X^{(m)} (Y) \, ( \Delta u(X) - \Delta u(Y))\,  d\mu(Y) &= \displaystyle\frac{2}{d} \, d^{-m}\Delta u(X) - \left(\frac{d+2}{d}\right)^{m}\,  \Delta_m u(X)
	\end{align*}

	\noindent  Then:

	\begin{align*}
	\Delta u(X) -\displaystyle  \frac{d}{2}\, (d+2)^{m} \Delta_m u(X) &= \frac{d}{2} \, d^{m} \int_{\mathfrak{SS}} \psi_X^{(m)} (Y) ( \Delta u(X) - \Delta u(Y)) \, d\mu(Y)
	\end{align*}

	\noindent Let us now consider the case where the Laplacian of~$u$ is H\"older continuous, with the exponent~$\alpha$ introduced above:

	$$\forall\, (X,Y)\,\in \mathfrak{SS}^2,\,\forall\,t\,\in\,[0,T] : \quad |\Delta u(t, X)-\Delta u(t, Y)| \leq C_\Delta  \, |X-Y|^\alpha$$
	\noindent where~$C_\Delta$ is a nonnegative real constant. This is the case, for instance, as soon as~$\Delta u(t,\cdot)$ is of finite energy, uniformly with respect to~$t$.\\

	\noindent One may note that, since~$\psi_X^{(m)}$ is nonnegative and since~$\displaystyle\int_{\mathfrak{SS}} \psi_X^{(m)}\,d\mu=\frac{2}{d}\,d^{-m}$:
	$$\displaystyle\frac{d}{2} d^{m} \int_{\mathfrak{SS}} \psi_X^{(m)} (Y) \, ( \Delta u(X) - \Delta u(Y))\, d\mu(Y)$$
	\noindent is a weighted mean value of~$\Delta u(X) - \Delta u(Y)$ over the support of~$\psi_X^{(m)}$, i.e. over the two~$m^{th}$-order cells containing $X$. The diameter of these cells is bounded by~$2^{-m}\,\delta$, which yields:

	\begin{align*}
	\left |\Delta u(X) - \displaystyle\frac{d}{2}(d+2)^{m}\,  \Delta_m u(X)\right | &\leq \frac{d}{2} \, d^{m} \int_{\mathfrak{SS}} \psi_X^{(m)} (Y) \,\left| \Delta u(X) - \Delta u(Y)\right| \, d\mu(Y)\\
	&\leq  C_\Delta\,\left(2^{-m}\,\delta\right)^{\alpha}\\
	&=   C_\Delta\,\delta^\alpha\,\left(\displaystyle\frac{1}{2}\right)^{m\,\alpha}
	\end{align*}

	\noindent  Finally:

	\begin{align*}
	\Delta u(X) &= \displaystyle \frac{d}{2}\, (d+2)^{m} \, \Delta_m u(X) + \mathcal{O}(2^{-m\,\alpha})
	\end{align*}

	\vskip 1cm
	\paragraph{Consistency}

	\begin{definition}
		The consistency error of the scheme is defined, for~$0\leq k \leq N-1$ and~$1\leq i \leq \mathcal{N}_m-d$, by:
		$$\varepsilon^m_{k,i}=\frac{u((k+1)\,h,X_i)-u(k\,h,X_i)}{h}-\frac{d}{2}\,(d+2)^m\,\Delta_m u(k\,h,X_i)$$
		The scheme is said to be \textbf{consistent} if the consistency error goes to zero when $h\rightarrow 0$ and $m \rightarrow+ \infty$, for some norm.
	\end{definition}

	\noindent Since~$u$ is a solution of the heat equation, the two estimates above yield:

	\[
	\varepsilon^m_{k,i} = \mathcal{O}(h) + \mathcal{O}(2^{-m \alpha}) \quad 0\leq k \leq N-1,\, 1\leq i \leq \mathcal{N}_m-d
	\]

	\noindent uniformly with respect to~$k$ and~$i$. One may check that

	$$ \displaystyle \lim_{h\rightarrow 0, \, m\rightarrow+\infty}\, \max_{k,i}\,|\varepsilon^m_{k,i}|=0$$

	\noindent The scheme is then consistent.

	\paragraph{Stability}

	\begin{definition}
		\noindent Let us recall that the \textbf{spectral norm}~$\rho$ is defined as the induced norm of the norm $\parallel \cdot \parallel_2$. It is given, for a square matrix~$A$, by:

		\[ \rho(A) =\sqrt{\lambda_{\max} \, \left( A^T \, A\right)}\]

		\noindent where $\lambda_{\max}$ stands for the spectral radius. For a symmetric matrix~$A$,~$\rho(A)$ is the largest absolute value of the eigenvalues of~$A$.
	\end{definition}

	\vskip 1cm

	\begin{proposition}

		\noindent Let us denote by~$\Phi$ the function such that:

		$$ \forall\, x \neq 0 :\quad  \Phi(x)=x\, (d+2-x).$$
		\noindent According to~\cite{Fukushima1992}, the Dirichlet eigenvalues~$\lambda_{m}$,~$m\,\in\,\N^\star$, of the graph Laplacians~$-\Delta_m$ (i.e. of the matrices~$L_m$) are related recursively:

		$$ \forall\, m\geq 2\, :\quad \lambda_{m-1}=\Phi(\lambda_{m}).$$

	\end{proposition}

	\vskip 1cm

	\noindent We deduce that, for any integer~$m\geq 2$:

	\[\lambda^{\pm}_{m}=\displaystyle\frac{(d+2)\pm \sqrt{(d+2)^2-4\, \lambda_{m-1}}}{2}\]

	\noindent Let us introduce the functions~$\phi^-$ and~$\phi^+$ such that, for any~$x$ in~\mbox{$ \left]-\infty,\displaystyle\frac{(d+2)^2}{4}\right]$} :

	$$
	\phi^-(x) =\displaystyle\frac{(d+2)- \sqrt{(d+2)^2-4\, x}}{2}\quad , \quad
	\phi^+(x)=\displaystyle\frac{(d+2) + \sqrt{(d+2)^2-4\, x}}{2}
	$$

	\noindent $\phi^+(0)=d+2$, $\phi^-\left(\frac{(d+2)^2}{4}\right)=\displaystyle\frac{d+2}{2}$, $\phi^-(0)=0$, and $\phi^+\left(\frac{(d+2)^2}{4}\right)=\displaystyle\frac{d+2}{2}$. In particular,~$\phi^-$ and~$\phi^+$ take their values in~$[0,d+2]$.\\

	\noindent The function~$\phi^-$ is increasing. Its fixed point is~$x^{-,\star}=0$.\\

	\noindent The function $\phi^+$ is non increasing. Its fixed point is~$x^{+,\star}=(d+2)-1$.\\
	\noindent One may also check that the following two maps are contractions in a neighbourhood of their fixed points, since:

	$$
	\left |\displaystyle \frac{d}{dx}\phi^-(0)\right | =\frac{1}{\sqrt{(d+2)^2}}
	=\frac{1}{d+2}<1$$

	\noindent and:

	$$
	\left |\displaystyle \frac{d}{d\, x}\phi^+\left ((d+2)-1\right )\right |=\frac{1}{\sqrt{(d+2)^2-4\, (d+2)+4}}
	=\displaystyle\frac{1}{d}<1.$$

	\noindent In~\cite{Shima},~T.~Shima gives the complete spectrum for $m\geq 1$. The Dirichlet eigenvalues of the graph Laplacian on $V_1$ are~$2$, with multiplicity $1$,~$d+2$, with multiplicity $d-1$, and~$2\,d$, with multiplicity~$\displaystyle\frac{d\,(d-3)}{2}$ (this last value does not appear in the case of the Gasket,~$d=3$).\\

	\noindent The complete Dirichlet spectrum, for $m\geq 2$, is generated by the maps $\phi^+$ and $\phi^-$, applied to the eigenvalues of the previous level, together with the initial values $2$, $d+2$ and $2\, d$.\\

	\noindent Since~$\phi^\pm$ take their values in~$[0,d+2]$, and since~$d+2\leq 2\,d$, one may finally conclude that, for any natural integer~$m$:

	\begin{align}
	\label{eqn:Spectralbound}
	0 < \lambda_m \leq 2\, d
	\end{align}

	\noindent This bound is optimal, since~$2\,d$ belongs to the spectrum for any~$m\geq 2$.

	\begin{remark}
		The eigenvalues~$\lambda$ of the Dirichlet Laplacian on~$\mathfrak{SS}$ are obtained as limits
		$$\lambda=\lim_{m\to+\infty}\frac{d}{2}\,(d+2)^m\,\lambda_m$$
		along sequences~$(\lambda_m)_m$ generated, from a certain rank, by~$\phi^-$ only. On the other hand, the eigenvalues of the matrix~$\tilde{\Delta}_m$ coming from the values~$2$,~$d+2$,~$2\,d$ created at level~$m$ are of order~$(d+2)^m$: these are the eigenvalues which constrain the time step of the explicit scheme.
	\end{remark}

	\vskip 1cm

	\begin{definition}
		\label{Stability}
		Let us set~$\delta x=2^{-m}$. The scheme is said to be:
		\begin{itemize}
			\item unconditionally stable if there exist a constant~$C \geq 1$ independent of $h$ and $m$ such that:

			$$\rho(A^{k})\leq C \quad \forall \, k\, \in \, \{1,\hdots,N\}$$

			\item conditionally stable if there exist three constants $\beta> 0$, $C_1 > 0$ and $C_2 \geq 1$, independent of $h$ and $m$, such that:

			$$ h \leq C_1 \, (\delta x)^{\beta} \Longrightarrow \rho(A^{k})\leq C_2  \quad \forall \, k\, \in \, \{1,\hdots,N\}$$

		\end{itemize}
	\end{definition}

	\vskip 1cm

	\begin{proposition}

		\noindent Let us denote by~$\gamma_i$,~\mbox{$i=1,\hdots,\mathcal{N}_m -d$}, the eigenvalues of the matrix~$A$. Then:

		\begin{equation}
		\label{eqn:ExplicitSpectrum}
		h\, (d+2)^{m}\leq \displaystyle\frac{2}{d^2} \Longrightarrow \forall\, i=1,\hdots,\mathcal{N}_m -d\,: \quad |\gamma_i|\leq 1.
		\end{equation}

		\noindent Since~$(d+2)^{-m}=(\delta x)^{D_W}$, where~$D_W=\displaystyle\frac{\ln(d+2)}{\ln 2}$ denotes the walk dimension of~$\mathfrak{SS}$, this condition writes:

		$$h\leq \frac{2}{d^2}\,(\delta x)^{D_W}$$

		\noindent The explicit scheme is thus conditionally stable, with~$\beta=D_W$,~$C_1=\displaystyle\frac{2}{d^2}$ and~$C_2=1$. We designate this result as the \textbf{Courant-Friedrichs-Lewy (CFL) stability condition}.\\

		\noindent Moreover, for~$m\geq 2$, this condition is also necessary: the exponent~$\beta=D_W$ is optimal.
	\end{proposition}

	\vskip 1cm

	\begin{proof}

		\noindent Let us recall our scheme writes, for any integer~$k$ belonging to~\mbox{$\left \lbrace 0, \hdots, N-1\right \rbrace $}:
		{
			$$\left ({\mathcal S}_{\mathcal H}\right) \quad \left \lbrace \begin{array}{cccc}
			\displaystyle \frac{u^m_h((k+1)\,h,X_i)-u^m_h(k\,h,X_i)}{h}&=&	\displaystyle \frac{d}{2}(d+2)^{m}  \,\displaystyle \sum_{Y \underset{m}\sim X_i}  \left ( u^m_h(k\,h,Y)- u^m_h(kh,X_i)\right) &  \quad \forall \,  1\leq i \leq \mathcal{N}_m-d \\
			u^m_h(k\,h,P_j)&=&0 &   \quad \forall \,  0\leq j \leq d-1  \\
			u^m_h(0,X_i)&=&g(X_i) & \quad 1\leq i \leq \mathcal{N}_m-d \\
			\end{array} \right.$$

			\noindent i.e., under matrix form:

			$$
			U^m_h(k) =
			\left(
			\begin{matrix}
			u^m_h(k\,h,X_{1})\\
			\vdots\\
			u^m_h(k\,h,X_{\mathcal{N}_m-d})\\
			\end{matrix}
			\right)
			\qquad \forall\,  k \, \in \left \lbrace 0, \hdots, N\right \rbrace.$$
		}
		\noindent It satisfies the recurrence relation:

		$$U^m_h(k+1) = A \, U^m_h(k) \qquad \forall\,  k \, \in \left \lbrace 0, \hdots, N-1\right \rbrace$$

		\noindent where:

		$$
		A=I_{\mathcal{N}_m -d}-h\,\tilde{\Delta}_{m}=I_{\mathcal{N}_m -d}-h\,\frac{d}{2}\,(d+2)^m\,L_m.
		$$

		\noindent One may use the recurrence to find:

		$$U^m_h(k)=A^{k}\, U^m_h(0) \qquad \forall\,  k \, \in \, \left \lbrace 1, \hdots, N\right \rbrace.$$

		\noindent The matrix~$A$ is symmetric. Its eigenvalues~$\gamma_i$,~\mbox{$i=1,\hdots,\mathcal{N}_m -d$}, are such that:

		$$\gamma_i = 1 - h\, \left( \frac{d}{2}(d+2)^{m} \right)\lambda_i$$

		\noindent where the~$\lambda_i$ denote the eigenvalues of~$L_m$. Due to the spectral bound~(\ref{eqn:Spectralbound}), one has, for any integer~$i$ belonging to~\mbox{$\left \lbrace1,\hdots,\mathcal{N}_m -d\right \rbrace$} :

		$$ 1 - h\, \displaystyle \frac{d}{2}(d+2)^{m}\, (2\, d) \leq\gamma_i< 1$$

		\noindent which leads to:

		$$h\, (d+2)^{m}\leq \displaystyle\frac{2}{d^2} \Longrightarrow |\gamma_i|\leq 1.$$

		\noindent Since~$A$ is symmetric, this yields~$\rho(A^k)=\rho(A)^k\leq 1$ for any~$k\,\in\,\{1,\hdots,N\}$.\\

		\noindent Conversely, for~$m\geq 2$, the value~$2\,d$ belongs to the spectrum of~$L_m$, so that the smallest eigenvalue of~$A$ is~$\gamma_{\min}=1-h\,d^2\,(d+2)^m$. If~$h\,(d+2)^m>\displaystyle\frac{2}{d^2}$, then~$\gamma_{\min}<-1$, and~$\rho(A^k)=|\gamma_{\min}|^k\to+\infty$: the scheme is unstable. The CFL condition is thus necessary and sufficient, and no exponent~$\beta<D_W$ can ensure stability.

	\end{proof}

	\vskip 1cm

	\begin{remark}[\textbf{Spectral, walk and self-similar dimensions}]
		\noindent Let us set~$\delta x=2^{-m}$, the diameter of the~$m^{th}$-order cells (up to the factor~$\delta$), and~$\delta t=h$. Let us recall that the self-similar (Hausdorff) dimension~$D_H$, the walk dimension~$D_W$ and the spectral dimension~$d_S$ of~$\mathfrak{SS}$ are given by (see~\cite{Kigami2001},~\cite{StrichartzLivre2006}):

		$$D_H=\frac{\ln d}{\ln 2}\quad , \quad D_W=\frac{\ln (d+2)}{\ln 2}\quad , \quad d_S=\frac{2\,\ln d}{\ln (d+2)}=\frac{2\,D_H}{D_W}$$

		\noindent These three dimensions appear naturally in the scheme:

		\begin{itemize}
			\item the measure of an~$m^{th}$-order cell is~$d^{-m}=(\delta x)^{D_H}$, and the number of unknowns satisfies~$\mathcal{N}_m-d\asymp d^m=(\delta x)^{-D_H}$;
			\item the normalization constant of the Laplacian is~$\displaystyle\frac{d}{2}\,(d+2)^m=\frac{d}{2}\,(\delta x)^{-D_W}$, and the largest eigenvalue of~$\tilde{\Delta}_m$ is~$d^2\,(d+2)^m=d^2\,(\delta x)^{-D_W}$;
			\item the CFL condition~(\ref{eqn:ExplicitSpectrum}) shows that the stability exponent is the walk dimension: $\beta=D_W$ replaces the exponent~$2$ of the classical CFL condition for the heat equation on a smooth domain (note that~$D_W>2$);
			\item since~$(d+2)^m=\left(d^m\right)^{\frac{2}{d_S}}$, the CFL condition can also be expressed in terms of the number of unknowns:
			$$\delta t \lesssim \left(\mathcal{N}_m\right)^{-\frac{2}{d_S}}$$
			which is the analogue of the classical condition~$\delta t\lesssim N^{-\frac{2}{n}}$ for~$N$ grid points in a domain of~$\R^n$: the dimension~$n$ is replaced by the spectral dimension~$d_S$. This is consistent with the Weyl asymptotics: the number of eigenvalues of~$\tilde{\Delta}_m$ less than~$\lambda$ behaves like~$\lambda^{\frac{d_S}{2}}$;
			\item the H\"{o}lder exponent satisfies~$\alpha=\displaystyle\frac{\ln\left(\frac{d+2}{d}\right)}{2\,\ln 2}=\frac{D_W-D_H}{2}$, so that the consistency error writes:
			$$\varepsilon^m_{k,i}=\mathcal{O}(\delta t)+\mathcal{O}\left((\delta x)^{\frac{D_W-D_H}{2}}\right)$$
		\end{itemize}
	\end{remark}

	\vskip 1cm

	\paragraph{Convergence}

	\begin{definition}
		\noindent
		\begin{itemize}
			\item The scheme is said to be convergent for the matrix norm $\| \cdot  \|$ if :
			{
				$$\displaystyle\lim_{h\rightarrow 0, \, m\rightarrow +\infty} \left  \| \left( u(kh,X_i)- u^m_h(k\,h,X_i) \right)_{0\leq k\leq N,\, 1\leq i \leq \mathcal{N}_m} \right  \|=0$$
			}
			\item The scheme is said to be conditionally convergent for the matrix norm $ \| \cdot  \|$ if there exist two real constants $\beta>0$ and $C>0$ such that, with~$\delta x=2^{-m}$ :

			$$\displaystyle \lim_{h\leq C\,(\delta x)^{\beta}, \, m\rightarrow +\infty} \left  \| \left( u(k\,h,X_i)- u^m_h(k\,h,X_i) \right)_{0\leq k\leq N,\, 1\leq i \leq \mathcal{N}_m}\right  \|=0$$

		\end{itemize}
	\end{definition}

	\vskip 1cm

	\begin{theorem}
		{	\noindent If the scheme is stable and consistent, then it is also convergent for the norm~$\| \cdot  \|_{2,\infty}$, such that:

			$$	\left \| \left( u^m_h(k\, h,X_i) \right)_{0\leq k\leq N, 1\leq i\leq \mathcal{N}_m} \right  \|_{2,\infty}=\displaystyle \max_{0\leq k\leq N} \left( d^{-m}\sum_{1\leq i\leq \mathcal{N}_m}  \left  |u^m_h(k\,h,X_i) \right |^2\right)^{\frac{1}{2}}$$

			\noindent In the case of the explicit scheme, under the CFL condition~(\ref{eqn:ExplicitSpectrum}), the scheme is conditionally convergent, with~$\beta=D_W$, and the error is of order~$\mathcal{O}(2^{-m\alpha})$.
		}
	\end{theorem}

	\vskip 1cm

	\begin{proof}

		\noindent Let us set:

		$$ w^k_i  = u(k\,h,X_i)-u^m_h(k\,h,X_i), \quad 0\leq k \leq N,\, 1\leq i\leq \mathcal{N}_m$$

		\noindent where the points of~$V_0$ are numbered~$X_{\mathcal{N}_m-d+1},\hdots,X_{\mathcal{N}_m}$. One may check that:

		\begin{align*}
		w^k_{\mathcal{N}_m-d+1}=\dots=w^k_{\mathcal{N}_m}&=0 \qquad 0\leq k \leq N\\
		w^0_i&=0 \qquad 1\leq i \leq \mathcal{N}_m -d
		\end{align*}

		\noindent Let us now introduce, for any integer~$k$ belonging to~\mbox{$\left \lbrace 0, \hdots,   N \right \rbrace$}:

		$$W^k=
		\left(
		\begin{matrix}
		w^k_1\\
		\vdots\\
		w^k_{\mathcal{N}_m-d}
		\end{matrix}
		\right)
		\quad , \quad  E^k=
		\left(
		\begin{matrix}
		\varepsilon^m_{k,1}\\
		\vdots\\
		\varepsilon^m_{k,\mathcal{N}_m-d}
		\end{matrix}
		\right)
		$$

		\noindent By definition of the consistency error, one has then~$W^0=0$, and, for any integer~$k$ belonging to~\mbox{$\left \lbrace 0, \hdots,   N-1\right \rbrace$}:

		$$
		W^{k+1} =A\, W^k+h\, E^k  $$

		\noindent One finds recursively, for any integer~$k$ belonging to~$\left \lbrace 1, \hdots,   N\right \rbrace$:

		$$
		W^{k}=A^k W^0+h\, \displaystyle\sum_{j=0}^{k-1} A^{k-1-j}\,  E^{j}=h\, \displaystyle\sum_{j=0}^{k-1} A^{k-1-j} \, E^{j}
		$$

		\noindent Since the matrix~$A$ is a symmetric one, the~\textbf{CFL} stability condition $h\,(d+2)^{m}\leq \displaystyle\frac{2}{d^2}$ yields~$\rho(A)\leq 1$, and, for any integer~$k$ belonging to~\mbox{$\left \lbrace  1,\hdots,N\right \rbrace $}:

		\begin{align*}
		|W^k|&\leq h\, \left(\sum_{j=0}^{k-1}\rho(A)^{k-1-j} \right)\, \left(\displaystyle\max_{0\leq j \leq N-1} | E^{j} | \right)\\
		&\leq h\, k\, \left(\max_{0\leq j \leq N-1} | E^{j} | \right)\\
		&\leq h\, N\, \left(\max_{0\leq j \leq N-1} | E^{j} | \right)\\
		&= T\, \left(\max_{0\leq j \leq N-1} \left(\sum_{i=1}^{\mathcal{N}_m-d} \, |\varepsilon^m_{j,i}|^2 \right)^{\frac{1}{2}} \right)
		\end{align*}

		\noindent One deduces then:

		\begin{align*}
		\max_{0\leq k\leq N} \left( d^{-m}\sum_{i=1}^{\mathcal{N}_m}   |w^k_i|^2\right)^{\frac{1}{2}}&=d^{-\frac{m}{2}}\max_{1\leq k\leq N}|\, W^k|\\
		&\leq \left(d^{-\frac{m}{2}}\right)\,  T\left(\max_{0\leq j \leq N-1} \left(\displaystyle\sum_{i=1}^{\mathcal{N}_m-d} |\varepsilon^m_{j,i}|^2 \right)^{\frac{1}{2}} \right)\\
		&\leq \left(d^{-\frac{m}{2}}\right) \, T\left((\mathcal{N}_m-d)^{\frac{1}{2}}\max_{0\leq j \leq N-1,\, 1\leq i \leq \mathcal{N}_m-d} |\varepsilon^m_{j,i}| \right)\\
		&=\sqrt{\left(d^{-m}\, \displaystyle\frac{d^{m+1}-d}{2}\right)} \, T\left(\displaystyle\max_{0\leq j \leq N-1,\, 1\leq i \leq \mathcal{N}_m-d} |\varepsilon^m_{j,i}| \right)\\
		&\leq \sqrt{\frac{d}{2}}\, T\left(\displaystyle\max_{0\leq j \leq N-1,\, 1\leq i \leq \mathcal{N}_m-d} |\varepsilon^m_{j,i}| \right)\\
		&=   \mathcal{O}(h) + \mathcal{O}(2^{-m \alpha})  \\
		&=   \mathcal{O}((d+2)^{-m}) + \mathcal{O}(2^{-m\alpha})  \\
		&= \mathcal{O}(2^{-m\alpha})
		\end{align*}

		\noindent since, under the CFL condition,~$h=\mathcal{O}((d+2)^{-m})$, and since~$\alpha<1$ yields~$(d+2)^{-m}\leq 2^{-m}\leq 2^{-m\alpha}$. The scheme is thus convergent.
	\end{proof}

	\vskip 1cm

	\begin{remark}
		\noindent One has to bear in mind that each~$m^{th}$-order cell has measure~$d^{-m}$, and that each point of~$V_m\setminus V_0$ belongs to exactly two~$m^{th}$-order cells. Thus, for any continuous function~$v$ on~$\mathfrak{SS}$ vanishing on~$V_0$:
		{
			$$\lim_{m\to+\infty}\frac{2}{d}\,d^{-m}\, \displaystyle \sum_{1\leq i\leq \mathcal{N}_m}   |v(X_i)|^2
			=\left \|  v \right \| ^2_{L^2(\mathfrak{SS},\mu)}$$

		}
		\noindent The norm~$\| \cdot  \|_{2}$ is thus, up to the factor~$\sqrt{d/2}$, a discrete analogue of the~$L^2$ norm.
	\end{remark}

	\vskip 1cm

		\subsubsection{The specific case of the implicit Euler Method}

		\noindent Let us consider the implicit Euler scheme, for any integer~$k$ belonging to~$\left \lbrace 1, \hdots, N \right \rbrace$, any point~$P_j$ of~$V_0$,~$0\leq j \leq d-1$, and any~$X$ in the set~\mbox{$V_m\setminus V_0$} :

		$$\left ({\mathcal S}_{\mathcal H}^{\rm imp}\right) \quad \left \lbrace \begin{array}{cccc}
		\displaystyle \frac{u^m_h(k\,h,X)-u^m_h((k-1)\,h,X)}{h}&=& \displaystyle \frac{d}{2}\,(d+2)^m \,\displaystyle \sum_{Y \underset{m}\sim X} \left( u^m_h(k\,h,Y)- u^m_h(kh,X)\right) &   \\
		u^m_h(k\,h,P_j)&=&0 &     \\
		u^m_h(0,X)&=&g(X) &
		\end{array} \right.$$

		\noindent where, as previously,~$\displaystyle\frac{d}{2}\,(d+2)^m=r^{-m}\left(\int_{\mathfrak{SS}} \psi^{(m)}_{X} d\mu\right)^{-1}$. Let us define the approximate equation as:

		$$
		u^m_h(k\,h,X) - h\, \frac{d}{2}\,(d+2)^m \, \sum_{Y \underset{m}\sim X} \left( u^m_h(k\, h,Y)- u^m_h(k\,h,X)\right)  = u^m_h((k-1)h,X) \quad
		\forall \, k\, \in\,\left \lbrace 1, \hdots, N \right \rbrace,  \, \forall \,X\,\in \,V_m\setminus V_0
		$$

		\noindent As before, we fix~$m\,\in \,\N$, and use the natural ordering~$X_1,\hdots,X_{\mathcal{N}_m-d}$ of the points of~$V_m\setminus V_0$ induced by their addresses~$X=f_w(P_i)$, where~\mbox{$w\,\in\,\{1,\dots,d\}^m$} denotes a word of length $m$, and where $P_i$,~$0 \leq i \leq d-1$, belongs to~$V_0$.\\

		\noindent We get, for any integer~$k$ belonging to~$ \left \lbrace 0, \hdots, N \right \rbrace$, the solution vector $U(k)$ as before:

		\begin{align*}
		U^m_h(k) &=\left(
		\begin{matrix}
		u^m_h(k\,h,X_1)\\
		\vdots\\
		u^m_h(k\,h,X_{\mathcal{N}_m-d})\\
		\end{matrix}
		\right)\\
		\end{align*}

		\noindent It satisfies the recurrence relation:

		$$\tilde{A} \, U^m_h(k) = U^m_h(k-1) \qquad \forall\,k\,\in\,\{1,\hdots,N\}$$

		\noindent where:

		$$
		\tilde{A}=I_{\mathcal{N}_m -d} + h\,\tilde{\Delta}_{m}=I_{\mathcal{N}_m -d} + h\,\frac{d}{2}\,(d+2)^m\,L_m
		$$

		\noindent  and where~$I_{\mathcal{N}_m -d}$ denotes the~$ {(\mathcal{N}_m -d)}\times  {(\mathcal{N}_m -d)}$ identity matrix, and $\tilde{\Delta}_{m}$ the~$ {(\mathcal{N}_m -d)}\times  {(\mathcal{N}_m -d)}$ normalized Laplacian matrix.\\

		\noindent The scheme is well defined. Indeed, the matrix~$\tilde{A}$ is symmetric; let us denote by~$\tilde{\gamma}_i$,~\mbox{$i=1,\hdots,\mathcal{N}_m -d$}, its eigenvalues. The spectral bound~(\ref{eqn:Spectralbound}) yields:

		$$1 < \tilde{\gamma}_i\leq 1+h\,d^2\,(d+2)^m$$

		\noindent The matrix~$\tilde{A}$ is thus invertible, and:

		\begin{align}
		\label{eqn:ImplicitSpectrum}
		\rho\left(\tilde{A}^{-1}\right)=\max_{1\leq i\leq \mathcal{N}_m -d}\,\left|\tilde{\gamma}_i^{-1}\right| < 1
		\end{align}

		\vskip 1cm

		\paragraph{Consistency, stability and convergence}

		\subparagraph{\emph{i}. The scheme error}

		\noindent Let~$u$ be the solution of the problem, and assume, as in the explicit case, that~$\displaystyle\frac{\partial^2 u}{\partial t^2}$ is bounded on~$[0,T]\times\mathfrak{SS}$. For all~$k$ in~$\left \lbrace 1, \hdots, N \right \rbrace$ :

		$$   \forall \,X\,\in \, {\mathfrak{SS}} \,:   \quad
		\displaystyle\frac{\partial u}{\partial t}(kh,X)=\displaystyle\frac{1}{h}\,\left( u(kh,X)-u((k-1)h,X)\right)+{\mathcal O}(h)
		$$

		\noindent On the other hand, if~$\Delta u$ is H\"{o}lder continuous of exponent~$\alpha$, uniformly with respect to the time variable, one has, for $X\in V_m\setminus V_0$:

		\begin{align*}
		\Delta u(kh,X) &= \frac{d}{2}(d+2)^{m} \Delta_m u(kh,X) + \mathcal{O}(2^{-m\alpha})
		\end{align*}

		\subparagraph{\emph{ii}. Consistency}

		\noindent The consistency error of the implicit Euler scheme is defined, for~$1\leq k \leq N$ and~$1\leq i \leq \mathcal{N}_m-d$, by:

		$$\varepsilon^m_{k,i}=\frac{u(k\,h,X_i)-u((k-1)\,h,X_i)}{h}-\frac{d}{2}\,(d+2)^m\,\Delta_m u(k\,h,X_i)$$

		\noindent Since~$u$ is a solution of the heat equation, one has:

		\[
		\varepsilon^m_{k,i} = \mathcal{O}(h) + \mathcal{O}(2^{-m\alpha}) \quad 1\leq k \leq N,\, 1\leq i \leq \mathcal{N}_m-d
		\]

		\noindent uniformly with respect to~$k$ and~$i$. We can check that

		\[ \lim_{h\rightarrow 0, m\rightarrow\infty}\,\max_{k,i}\, |\varepsilon^m_{k,i}|=0\]

		\noindent The scheme is then consistent.

		\paragraph{Stability}

		\noindent Let us recall that our scheme writes:
		\footnotesize

		$$\left ({\mathcal S}_{\mathcal H}^{\rm imp}\right) \quad \left \lbrace \begin{array}{cccc}
		\displaystyle \frac{u^m_h(k\,h,X_i)-u^m_h((k-1)\,h,X_i)}{h}&=& \displaystyle \frac{d}{2}(d+2)^{m}  \,\displaystyle \sum_{Y \underset{m}\sim X_i}  \left(u^m_h(k\,h,Y)- u^m_h(kh,X_i)\right) &  \quad \forall\,  1\leq k\leq N,\, 1\leq i \leq \mathcal{N}_m-d \\
		u^m_h(k\,h,P_j)&=&0 &   \quad \forall \, 1\leq k\leq N,\,0\leq j\leq d-1  \\
		u^m_h(0,X_i)&=&g(X_i) & \quad 1\leq i \leq \mathcal{N}_m-d \\
		\end{array} \right.$$
		\normalsize

		\noindent i.e., under matrix form :

		$$U^m_h(k)=\tilde{A}^{-1} \, U^m_h(k-1)=\tilde{A}^{-k} \, U^m_h(0) \qquad \forall\,k\,\in\,\{1,\hdots,N\}$$

		\noindent Since~$\tilde{A}^{-1}$ is symmetric, the spectral bound~(\ref{eqn:ImplicitSpectrum}) yields, for any~$h>0$ and any~$m\,\in\,\N$:

		$$\rho\left(\tilde{A}^{-k}\right)=\rho\left(\tilde{A}^{-1}\right)^k\leq 1 \quad \forall\,k\,\in\,\{1,\hdots,N\}$$

		\noindent and thus, denoting by~$|\cdot|$ the Euclidean norm:

		$$\left|U^m_h(k)\right|\leq \left|U^m_h(0)\right| \quad \forall\,k\,\in\,\{1,\hdots,N\}$$

		\noindent This enables us to conclude that, in the sense of definition~\ref{Stability}, the implicit scheme is unconditionally stable, with~$C=1$: contrary to the explicit case, no condition linking~$h$ and~$\delta x=2^{-m}$ is required.\\

		\noindent One may also note that~$\tilde{A}$ is a strictly diagonally dominant M-matrix. Thus,~$\tilde{A}^{-1}$ has nonnegative coefficients, and:
		$$ \parallel \tilde{A}^{-1} \parallel_{\infty} \leq 1\quad \text{and thus} \quad
		\parallel \tilde{A}^{-k} \parallel_{\infty} \leq 1 \quad \forall\,k\,\in\,\{1,\hdots,N\}$$

		\noindent i.e.~$\parallel U^m_h(k)\parallel_{\infty} \leq \parallel U^m_h(0) \parallel_{\infty}$: the implicit scheme also satisfies the discrete maximum principle.

		\vskip 1cm

		\subparagraph{\emph{iii}. Convergence}

		\begin{theorem}
			\noindent The implicit Euler scheme is convergent for the norm $\parallel . \parallel_{2,\infty}$, without any condition linking~$h$ and~$m$:

			$$\left\| \left( u(kh,X_i)- u^m_h(k\,h,X_i) \right)\right\|_{2,\infty}=\mathcal{O}(h)+\mathcal{O}(2^{-m\alpha})$$

		\end{theorem}

		\vskip 0.5cm

		\begin{proof}

			\noindent Let:

			$$ w^k_i = u(kh,X_i)-u^m_h(kh,X_i), \quad 0\leq k \leq N,\, 1\leq i\leq \mathcal{N}_m$$

			\noindent One may check that:

			\begin{align*}
			w^k_{\mathcal{N}_m-d+1}=\dots=w^k_{\mathcal{N}_m}&=0 \qquad 0\leq k \leq N\\
			w^0_i&=0 \qquad 1\leq i \leq \mathcal{N}_m -d
			\end{align*}

			\noindent We set:

			$$ W^k=
			\left(
			\begin{matrix}
			w^k_1\\
			\vdots\\
			w^k_{\mathcal{N}_m-d}
			\end{matrix}
			\right)
			\quad , \quad
			E^k=
			\left(
			\begin{matrix}
			\varepsilon^m_{k,1}\\
			\vdots\\
			\varepsilon^m_{k,\mathcal{N}_m-d}
			\end{matrix}
			\right)
			$$

			\noindent By definition of the consistency error, the exact solution satisfies~$\tilde{A}\,\left(u(kh,X_i)\right)_{i}=\left(u((k-1)h,X_i)\right)_{i}+h\,E^k$. Thus,~$W^0=0$, and,  for~$1\leq k \leq N $:

			\begin{align*}
			W^{k}&=\tilde{A}^{-1} \left(W^{k-1}+h\, E^k\right)
			\end{align*}

			\noindent We find, by induction, for~$1\leq k \leq N $:

			\begin{align*}
			W^{k}&=\tilde{A}^{-k} W^0+h\, \displaystyle \sum_{j=1}^{k} \tilde{A}^{-(k-j+1)} E^{j}  \\
			&=h\, \displaystyle \sum_{j=1}^{k} \tilde{A}^{-(k-j+1)} E^{j}
			\end{align*}

			\noindent Due to the stability of the scheme,~$\rho(\tilde{A}^{-1})\leq 1$, and we have, for~$k=1,\hdots,N$:

			\begin{align*}
			|W^k|&\leq h\,\left(\displaystyle \sum_{j=1}^{k}\rho(\tilde{A}^{-1})^{k-j+1} \right)\left(\max_{1\leq j \leq N} | E^{j} | \right)\\
			&\leq h\, k\,\left(\max_{1\leq j \leq N} | E^{j} | \right)\\
			&\leq h\,N\,\left(\max_{1\leq j \leq N} | E^{j} | \right)\\
			&= T\,\left(\max_{1\leq j \leq N} \left(\sum_{i=1}^{\mathcal{N}_m-d} |\varepsilon^m_{j,i}|^2 \right)^{1/2} \right)
			\end{align*}

			\noindent One deduces then:

			\begin{align*}
			\max_{0\leq k\leq N} \left( d^{-m}\sum_{i=1}^{\mathcal{N}_m}   |w^k_i|^2\right)^{\frac{1}{2}}&=d^{-\frac{m}{2}}\max_{1\leq k\leq N}|W^k|\\
			&\leq \left(d^{-\frac{m}{2}}\right) T\left(\max_{1\leq j \leq N} \left(\sum_{i=1}^{\mathcal{N}_m-d} |\varepsilon^m_{j,i}|^2 \right)^{1/2} \right)\\
			&\leq \left(d^{-\frac{m}{2}}\right) T\left((\mathcal{N}_m-d)^{\frac{1}{2}}\max_{1\leq j \leq N,\, 1\leq i \leq \mathcal{N}_m-d} |\varepsilon^m_{j,i}| \right)\\
			&=\sqrt{\left(d^{-m}\frac{d^{m+1}-d}{2}\right)} T\left(\max_{1\leq j \leq N,\, 1\leq i \leq \mathcal{N}_m-d} |\varepsilon^m_{j,i}| \right)\\
			&\leq \sqrt{\frac{d}{2}}\, T\left(\max_{1\leq j \leq N,\, 1\leq i \leq \mathcal{N}_m-d} |\varepsilon^m_{j,i}| \right)\\
			&= \mathcal{O}(h) + \mathcal{O}(2^{-m\alpha})
			\end{align*}

			\noindent which goes to zero when~$h\to 0$ and~$m\to+\infty$, independently. The scheme is thus convergent.
		\end{proof}

		\vskip 1cm

		\begin{remark}
			\noindent With the notations of the remark on dimensions,~$\delta x=2^{-m}$, the error of the implicit scheme writes:

			$$\mathcal{O}(\delta t)+\mathcal{O}\left((\delta x)^{\frac{D_W-D_H}{2}}\right)$$

			\noindent The two contributions are balanced for~$\delta t\asymp (\delta x)^{\frac{D_W-D_H}{2}}$, which is much larger than the time step~$\delta t\lesssim (\delta x)^{D_W}$ imposed by the stability of the explicit scheme. For a given accuracy, the implicit scheme thus requires far fewer time steps, at the cost of the resolution of a sparse, symmetric positive definite linear system at each step.
		\end{remark}
	
	\vskip 1cm
	
\subsubsection{Numerical results - Gasket and Tetrahedron}

	\paragraph{Recursive construction of the matrix related to the sequence of graph Laplacians}$\,$\\

	In the sequel, we describe our recursive algorithm used to construct the matrix related to the sequence of graph Laplacians, in the case of Sierpi\'{n}ski Gasket and Tetrahedron. The matrix~$A_m$ obtained below is the matrix of~$-\Delta_m$ on the whole set~$V_m$; the matrix~$L_m$ of the scheme is then obtained by deleting the rows and the columns related to the points of~$V_0$.\\

	\subparagraph{\emph{i}. The Sierpi\'{n}ski Gasket.}

	\begin{center}
		\includegraphics[scale=0.5]{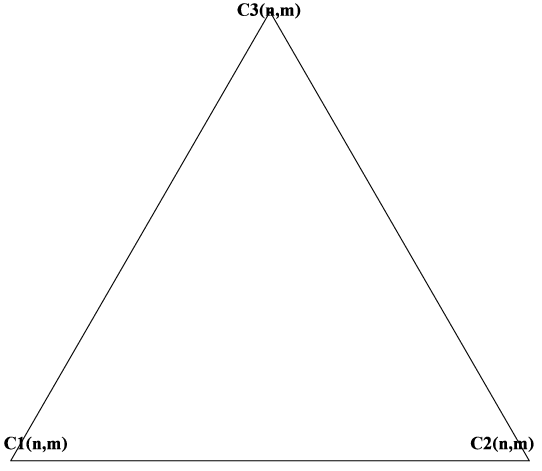}
		\captionof{figure}{$m^{th}$-order cell of the Sierpi\'{n}ski Gasket.}
		\label{figSGcell}
	\end{center}

	\noindent One may note, first, that, given a strictly positive integer~$m$,  a~$m^{th}$-order triangle has three corners, that we will denote by~$C1$,~$C2$ and~$C3$ ; the $(m+1)^{th}$-order triangle is then constructed by connecting three copies~$T(n)$,~$n=1,\,2,\, 3$, of the~$m^{th}$-order triangle.\\

	\noindent The initial triangle is labeled such that $C1 \sim 1$, $C2 \sim 2$ and $C3 \sim 3$~(see figure~\ref{figSGcell}).

	\begin{figure}[!htb]

		\minipage{0.32\textwidth}
		\includegraphics[width=\linewidth]{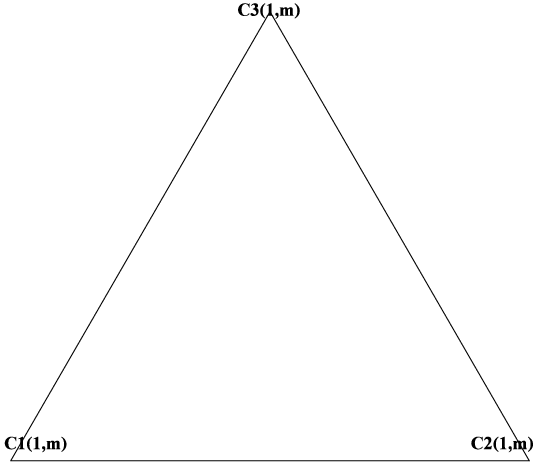}
		\caption{The first copy $T(1)$}
		\label{figSG1}
		\endminipage\hfill
		\minipage{0.32\textwidth}
		\includegraphics[width=\linewidth]{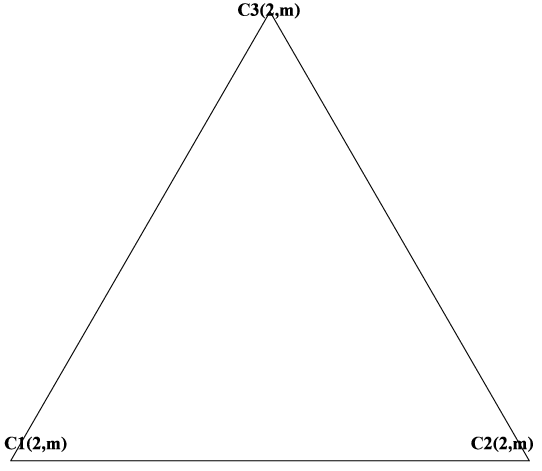}
		\caption{The second copy $T(2)$}
		\label{figSG2}
		\endminipage\hfill
		\begin{center}
			\minipage{0.32\textwidth}%
			\includegraphics[width=\linewidth]{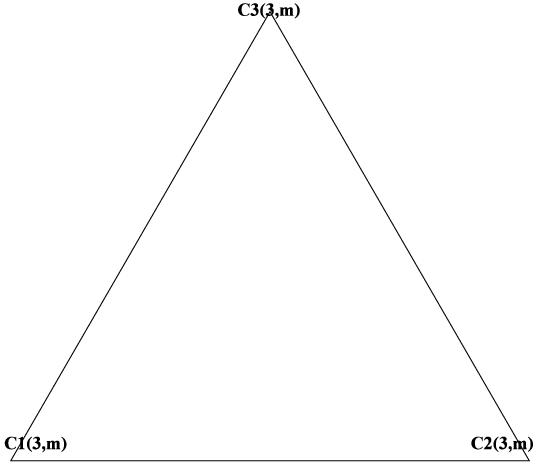}
			\caption{The third copy $T(3)$}
			\label{figSG3}
			\endminipage
		\end{center}
	\end{figure}

	\noindent The fusion is done by connecting $C2(1,m) \sim C1(2,m)$, $C3(1,m) \sim C1(3,m)$, and $C3(2,m) \sim C2(3,m)$~(see figures~\ref{figSG1},~\ref{figSG2},~\ref{figSG3}).\\

	\noindent The label of the corner vertex~$Ci$ of the copy~$T(n)$ can be obtained by means of the following recursive sequence, for any strictly positive integer~$m$:

	\begin{align*}
	C1(n,m)&=1+(n-1)\,\mathcal{N}_{m-1}\\
	C2(n,m)&=I2(m)+(n-1)\,\mathcal{N}_{m-1}\\
	C3(n,m)&=n \,\mathcal{N}_{m-1}\\
	\end{align*}

	\noindent where~$I2(m)$ denotes the label of the corner~$C2$ in the~$(m-1)^{th}$-order triangle, and where:

	\begin{align*}
	\mathcal{N}_{-1}&=3 \quad , \quad
	I2(0)=0\\
	I2(m)&=I2(m-1)+\mathcal{N}_{m-2}-1.\\
	\end{align*}

	\noindent The value~$\mathcal{N}_{-1}=3$ is a convention, which ensures that~$I2(1)=2$.
	\newpage
	\begin{enumerate}

		\item One may start with the initial triangle with the set of vertices~$V_0$. The corresponding matrix is given by:

		$$A_0 = \left(
		\begin{matrix}
		2 & -1 & -1 \\
		-1 & 2 & -1 \\
		-1 & -1 & 2 \\
		\end{matrix}
		\right)$$

		\item If $m=0$, the Laplacian matrix is $A_0$, else, $A_m$ is constructed recursively from three copies of the Laplacian matrices $A_{m-1}$ of the graph $V_{m-1}$. First, we build, for any strictly positive integer~$m$, the block diagonal matrix:

		$$B_m = \left(
		\begin{matrix}
		A_{m-1} & 0 & 0 \\
		0 & A_{m-1} & 0 \\
		0 & 0 & A_{m-1} \\
		\end{matrix}
		\right)$$

		\item One may then introduce, for any strictly positive integer~$m$, the connection matrix as in~\cite{UtaFreiberg2004}:

		$$C_m = \left(
		\begin{matrix}
		C2(1,m) & C3(1,m) & C3(2,m) \\
		C1(2,m) & C1(3,m) & C2(3,m) \\
		\end{matrix}
		\right)$$

		{	\item For each column~$j$ of~$C_m$, one has then to add the row (resp. the column)~$C_m(2,j)$ of~$B_m$ to the row (resp. the column)~$C_m(1,j)$. One then deletes the rows and the columns~$C_m(2,j)$, for all~$j$.}

	\end{enumerate}

	\vskip 1cm

	\subparagraph{\emph{ii}. The Sierpi\'{n}ski Tetrahedron.}

	\vskip 1cm

	\begin{center}
		\includegraphics[scale=0.7]{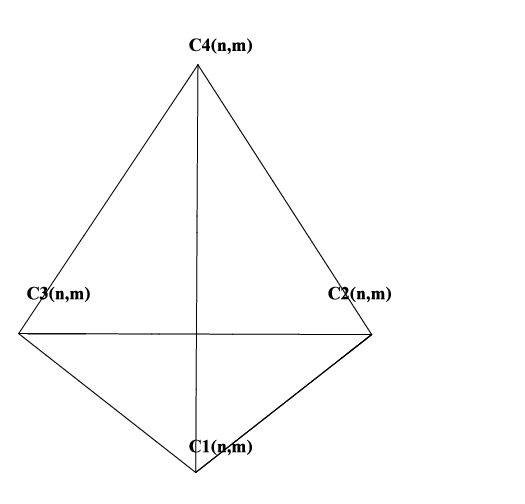}
		\captionof{figure}{$m^{th}$-order cell of the Sierpi\'{n}ski Tetrahedron.}
		\label{figSTcell}
	\end{center}

\begin{figure}[!htb]

		\minipage{0.32\textwidth}
		\includegraphics[width=\linewidth]{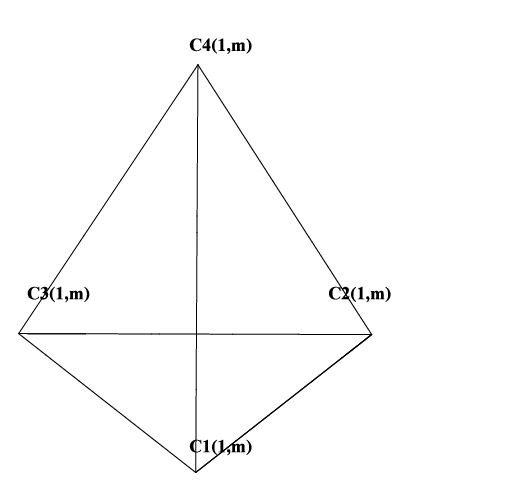}
		\caption{The first copy $T(1)$.}
		\label{figST1}
		\endminipage\hfill
		\minipage{0.32\textwidth}
		\includegraphics[width=\linewidth]{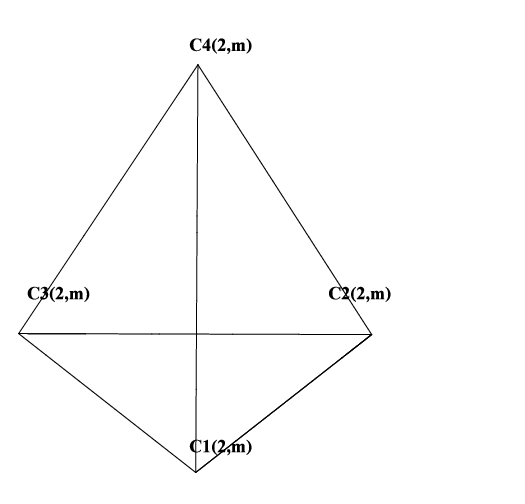}
		\caption{The second copy $T(2)$.}
		\label{figST2}
		\endminipage\hfill
		\minipage{0.32\textwidth}
		\includegraphics[width=\linewidth]{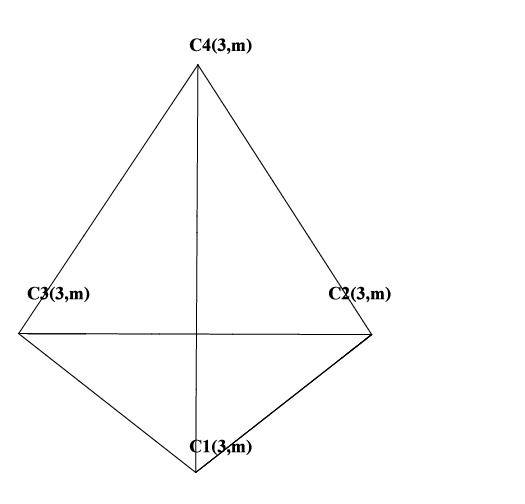}
		\caption{The third copy $T(3)$.}
		\label{figST3}
		\endminipage\hfill

		\begin{center}
			\minipage{0.32\textwidth}%
			\includegraphics[width=\linewidth]{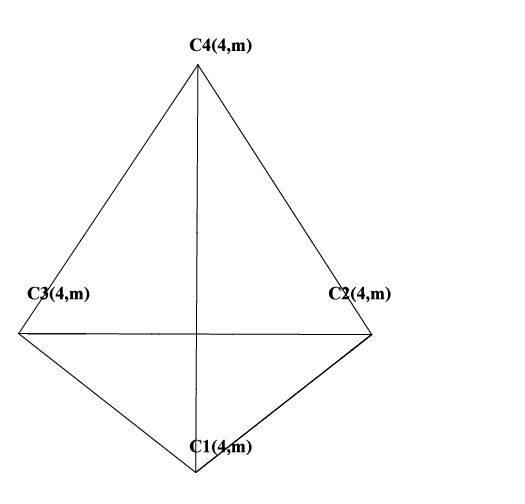}
			\caption{The fourth copy $T(4)$.}
			\label{figST4}
			\endminipage
		\end{center}
	\end{figure}

	One may note, first, that, given a strictly positive integer~$m$,  a~$m^{th}$-order tetrahedron has four corners~$C1$,~$C2$,~$C3$ and~$C4$~(see figure~\ref{figSTcell}), and that the $(m+1)^{th}$-order tetrahedron is constructed by connecting four copies~$T(n)$,~$n=1,2,3,4$, of the~$m^{th}$-order tetrahedron~(see figures~\ref{figST1} to~\ref{figST4}).\\

	As in the case of the triangle, the initial tetrahedron is labeled such that $C1 \sim 1$, $C2 \sim 2$, $C3 \sim 3$ and $C4 \sim 4$.

	\noindent The fusion is done by connecting $C2(1,m) \sim C1(2,m)$, $C3(1,m) \sim C1(3,m)$, $C4(1,m) \sim C1(4,m)$, $C3(2,m) \sim C2(3,m)$, $C4(2,m) \sim C2(4,m)$, $C4(3,m) \sim C3(4,m)$.\\

	\newpage
	\noindent The labels of the corners can be obtained by means of the following recursive sequence, for any strictly positive integer~$m$:

	\begin{align*}
	C1(n,m)&=1+(n-1)\,\mathcal{N}_{m-1}\\
	C2(n,m)&=I2(m)+(n-1)\,\mathcal{N}_{m-1}\\
	C3(n,m)&=I3(m)+(n-1) \,\mathcal{N}_{m-1}\\
	C4(n,m)&=n \,\mathcal{N}_{m-1}\\
	\end{align*}

	\noindent where:
	{
		\begin{align*}
		\mathcal{N}_{-1}&=3\\
		I2(0)&=0\\
		I2(m)&=I2(m-1)+\mathcal{N}_{m-2}-1\\
		I3(1)&=3\\
		I3(m)&=I3(m-1)+2\times\mathcal{N}_{m-2}-3\\
		\end{align*}
	}
	\noindent As before, the value~$\mathcal{N}_{-1}=3$ is a convention, which ensures that~$I2(1)=2$.
	\newpage
	\begin{enumerate}

		\item One starts with initial tetrahedron with the set of vertices $V_0$. The corresponding matrix is given by:

		$$A_0 = \left(
		\begin{matrix}
		3 & -1 & -1 & -1\\
		-1 & 3 & -1 & -1\\
		-1 & -1 & 3 & -1 \\
		-1 & -1 & -1 & 3 \\
		\end{matrix}
		\right)$$
		\vskip 0.5cm

		\item If $m=0$ the Laplacian matrix is $A_0$, else, for any strictly positive integer~$m$, $A_m$ is constructed recursively from four copies of the Laplacian matrices $A_{m-1}$ of the graph $V_{m-1}$. Thus, we build the block diagonal matrix:

		$$B_m = \left(
		\begin{matrix}
		A_{m-1} & 0 & 0 & 0 \\
		0 & A_{m-1} & 0 & 0 \\
		0 & 0 & A_{m-1} & 0 \\
		0 & 0 & 0 & A_{m-1} \\
		\end{matrix}
		\right)$$

		\vskip 0.5cm
		\item We then write the connection matrix:
		{
			$$C_m = \left(
			\begin{matrix}
			C2(1,m) & C3(1,m) & C3(2,m) & C4(1,m) & C4(2,m) & C4(3,m)\\
			C1(2,m) & C1(3,m) & C2(3,m) & C1(4,m) & C2(4,m) & C3(4,m)\\
			\end{matrix}
			\right)$$

			\vskip 0.5cm

			\item For each column~$j$ of~$C_m$, one then has to add the row (resp. the column)~$C_m(2,j)$ of~$B_m$ to the row (resp. the column)~$C_m(1,j)$. One then deletes the rows and the columns~$C_m(2,j)$, for all~$j$.}

	\end{enumerate}

	\vskip 1cm

	\newpage
	\paragraph{Numerical results}

	\noindent Our heat transfer simulation consists in a propagation scenario, where the initial condition is a harmonic spline~$g$, the support of which consists of the~$m^{th}$-order cells containing a vertex~$x$ of~$V_m\setminus V_0$, such that it takes the value~$1$ on~$x$, and~$0$ on the other vertices of~$V_m$. In both cases, the parameters satisfy the CFL condition~(\ref{eqn:ExplicitSpectrum}), i.e.~$h\leq \displaystyle\frac{2}{d^2}\,(\delta x)^{D_W}$, with~$\delta x=2^{-m}$.\\

	\noindent The colour scale represents the value of the temperature, from red (high values) to blue (low values).\\

	\subparagraph{\emph{i}. The Sierpi\'{n}ski Gasket}$\,$\\

	\noindent In the sequel (see figures~\ref{figSGk0} to~\ref{figSGk1000}), we present the numerical results for~$m=6$, $T=1$ and $N=2\times 10^5$, i.e.~$h=5\times 10^{-6}$, while the CFL condition requires~$h\leq \displaystyle\frac{2}{9\times 5^6}\approx 1.42\times 10^{-5}$.

	\begin{center}
		\includegraphics[scale=1.5]{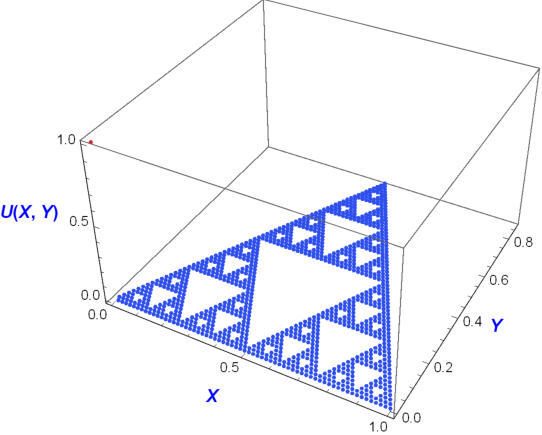}
		\captionof{figure}{The graph of the approximate solution of the heat equation for $k=0$.}
		\label{figSGk0}
	\end{center}

	\begin{center}
		\includegraphics[scale=1.5]{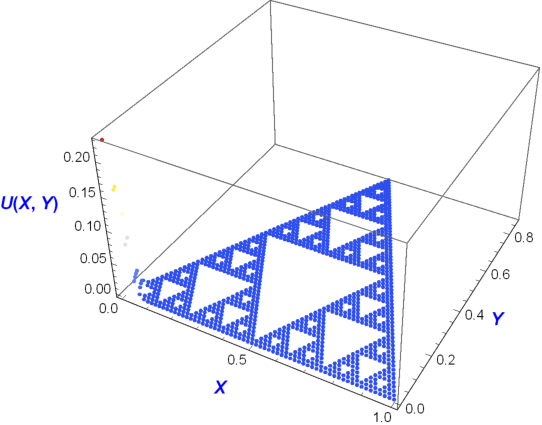}
		\captionof{figure}{The graph of the approximate solution of the heat equation for $k=10$.}
		\label{figSGk10}
	\end{center}

	\begin{center}
		\includegraphics[scale=1.5]{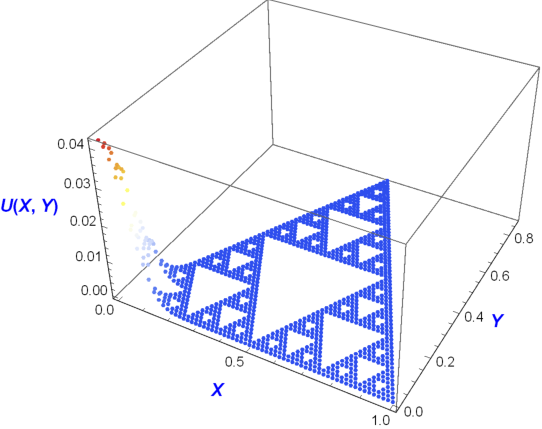}
		\captionof{figure}{The graph of the approximate solution of the heat equation for $k=100$.}
		\label{figSGk100}
	\end{center}

	\begin{center}
		\includegraphics[scale=1.5]{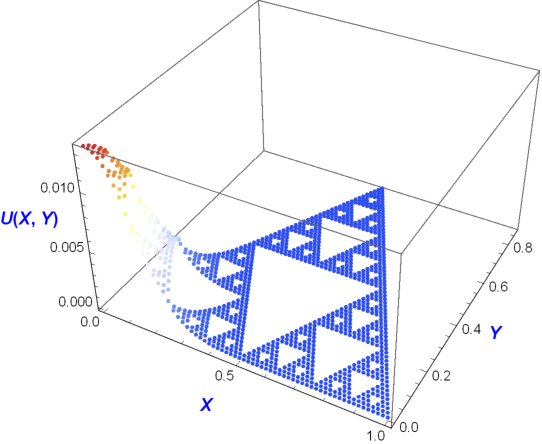}
		\captionof{figure}{The graph of the approximate solution of the heat equation for $k=500$.}
		\label{figSGk500}
	\end{center}

	\begin{center}
		\includegraphics[scale=1.5]{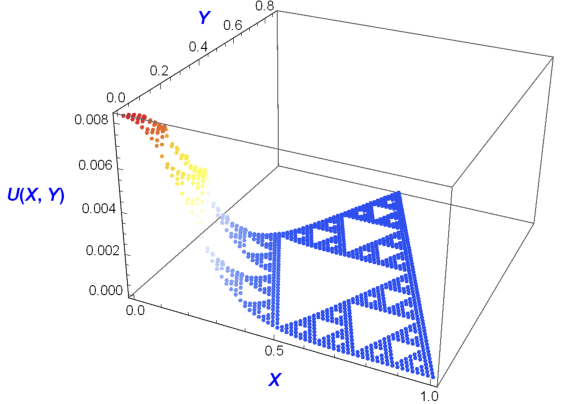}
		\captionof{figure}{The graph of the approximate solution of the heat equation for $k=1000$.}
		\label{figSGk1000}
	\end{center}

	\newpage
	\subparagraph{\emph{ii}. The Sierpi\'{n}ski Tetrahedron}$\,$\\

		In the sequel (see figures~\ref{figSTk0} to~\ref{figSTk500}), we present the numerical results for~$m=5$, $T=1$ and $N=10^5$, i.e.~$h=10^{-5}$, while the CFL condition requires~$h\leq \displaystyle\frac{2}{16\times 6^5}\approx 1.61\times 10^{-5}$.  \\

	\begin{center}
		\includegraphics[scale=1]{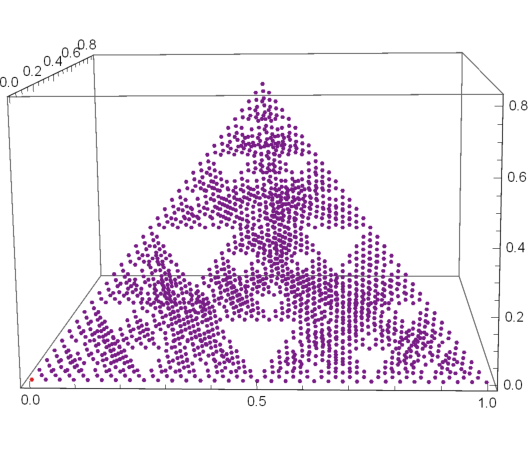}
		\captionof{figure}{The graph of the approximate solution of the heat equation for $k=0$.}
		\label{figSTk0}
	\end{center}

	\begin{center}
		\includegraphics[scale=1]{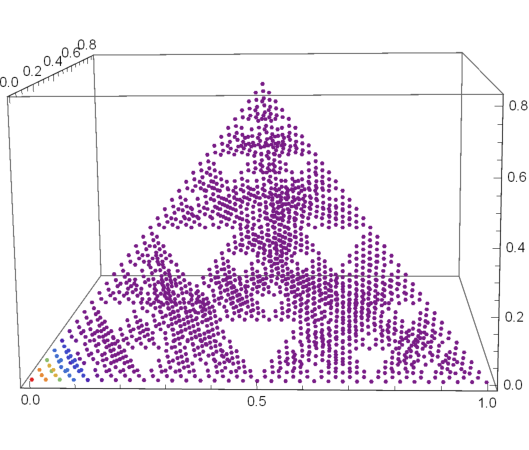}
		\captionof{figure}{The graph of the approximate solution of the heat equation for $k=10$.}
		\label{figSTk10}
	\end{center}

	\begin{center}
		\includegraphics[scale=1]{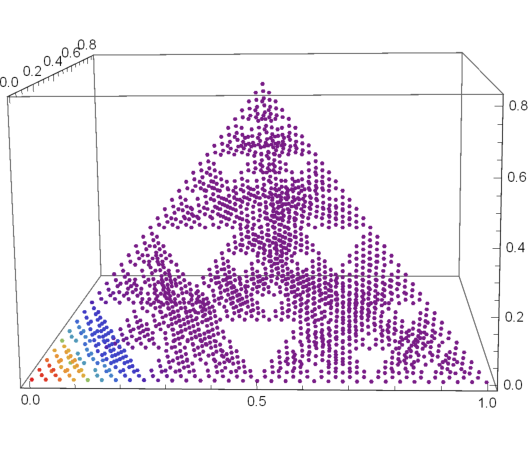}
		\captionof{figure}{The graph of the approximate solution of the heat equation for $k=50$.}
		\label{figSTk50}
	\end{center}

	\begin{center}
		\includegraphics[scale=1]{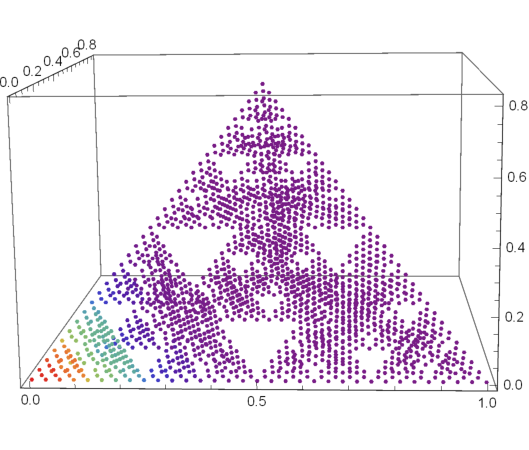}
		\captionof{figure}{The graph of the approximate solution of the heat equation for $k=100$.}
		\label{figSTk100}
	\end{center}

	\begin{center}
		\includegraphics[scale=1]{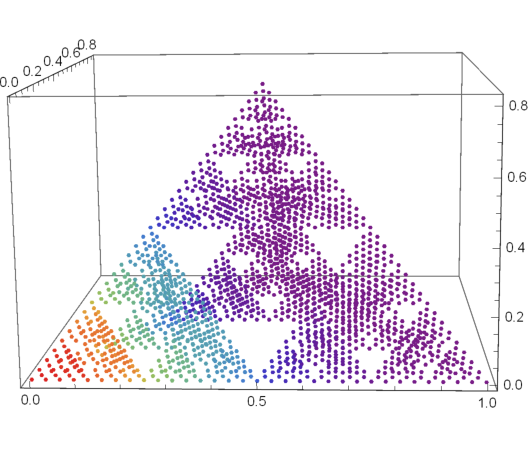}
		\captionof{figure}{The graph of the approximate solution of the heat equation for $k=500$.}
		\label{figSTk500}
	\end{center}

	\vskip 1cm

		An interesting feature in our work is that, contrary to existing ones, we do not rely on heat kernel estimates. Using a direct method has thus enabled us to discuss the choices of parameters as the integer~$m$, the step~$h$, and the convergence.\\

		As expected, the explicit numerical scheme is unstable and diverges unless one respects the stability condition between~$h$ and~$m$. Also, the propagation process evolves with time, directed from hot regions, towards cold ones.\\

		In order to go further, we have also studied the evolution in time of the temperature of a point $x_0$. Figures~\ref{figTemp} and~\ref{figLogTemp} respectively display the graph and log-log graph of the temperature as a function of time.\\


		\begin{center}
			\includegraphics[scale=1]{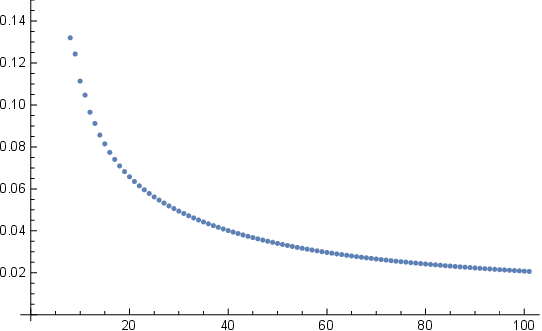}
			\captionof{figure}{The graph of the temperature as a function of time.}
			\label{figTemp}
		\end{center}

		\begin{center}
			\includegraphics[scale=1]{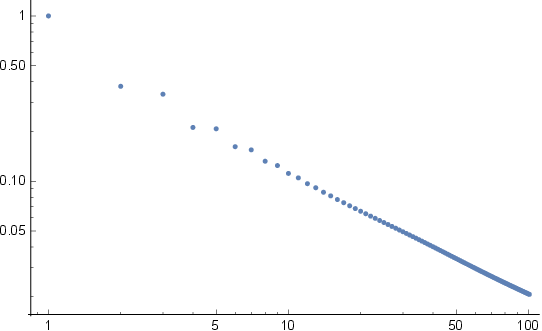}
			\captionof{figure}{The log-log graph of the temperature as a function of time.}
			\label{figLogTemp}
		\end{center}

		We find, numerically, that, on the time interval considered, the temperature follows the law:

		$$\ln \left  (u(t,x_0)\right ) \approx -1.01275 - 1.44972 \, \ln t$$

		\noindent i.e. a power law decay of the form:

		$$u(t,x_0)\approx C\, t^{-\kappa}\quad , \quad \kappa\approx 1.45$$

		\noindent where~$C$ is a strictly positive real constant. This result holds for different values of~$m$.\\

		\noindent Let us recall that, for small times, the heat kernel~$p_t$ on the Sierpi\'{n}ski Gasket satisfies the on-diagonal estimate~$p_t(x,x)\asymp t^{-\frac{d_S}{2}}$, where~$d_S=\displaystyle\frac{2\,\ln 3}{\ln 5}\approx 1.365$ is the spectral dimension. The exponent~$\kappa$ observed here differs from~$\displaystyle\frac{d_S}{2}\approx 0.68$. This is not contradictory: our initial condition is not a point mass, in our simulations, the heat source lies close to a boundary point of~$V_0$, where the Dirichlet condition accelerates the decay, and, for large times, the solution decays exponentially, like~$e^{-\lambda_1 t}$, where~$\lambda_1$ denotes the first Dirichlet eigenvalue. The observed power law thus only holds on an intermediate range of times. A finer study of the short-time regime, which should allow one to recover the spectral dimension numerically, will be the subject of further work.
		

		
		


		
		



	\bibliographystyle{amsplain}
	\bibliography{BibliographieClaire}

\providecommand{\bysame}{\leavevmode\hbox to3em{\hrulefill}\thinspace}
\providecommand{\MR}{\relax\ifhmode\unskip\space\fi MR }
\providecommand{\MRhref}[2]{%
  \href{http://www.ams.org/mathscinet-getitem?mr=#1}{#2}
}
\providecommand{\href}[2]{#2}
\begin{thebibliography}{10}

\bibitem{DalrympleFDM}
K.~Dalrymple, R.~S. Strichartz, and J.~P. Vinson, \emph{Fractal differential
  equations on the {S}ierpi\'{n}ski {G}asket}, The Journal of Fourier Analysis
  and Applications \textbf{5} (1999), no.~2/3, 203--284.

\bibitem{UtaFreiberg2004}
U.~R. Freiberg and M.~R. Lancia, \emph{Energy form on a closed fractal curve},
  Journal for Analysis and its Applications \textbf{23} (2004), no.~1,
  115--137.

\bibitem{Fukushima1992}
M.~Fukushima and T.~Shima, \emph{On a spectral analysis for the
  {S}ierpi\'{n}ski gasket}, Potential Anal \textbf{1} (1992), 1--3.

\bibitem{GibbonsFEM}
M.~Gibbons, A.~Raj, and R.~S. Strichartz, \emph{The {F}inite {E}lement {M}ethod
  on the {S}ierpi\'{n}ski gasket}, Constructive Approximation \textbf{17}
  (2001), no.~4, 561--588.

\bibitem{Hutchinson1981}
J.~E. Hutchinson, \emph{Fractals and self similarity}, Indiana University
  Mathematics Journal \textbf{30} (1981), 713--747.

\bibitem{Kigami1989}
J.~Kigami, \emph{A harmonic calculus on the {S}ierpi\' nski spaces}, Japan J.
  Appl. Math. \textbf{8} (1989), 259--290.

\bibitem{Kigami1993}
\bysame, \emph{Harmonic calculus on p.c.f. self-similar sets}, Trans. Amer.
  Math. Soc. \textbf{335} (1993), 721--755.

\bibitem{Kigami2001}
\bysame, \emph{Analysis on fractals}, Cambridge University Press, 2001.

\bibitem{Kigami2003}
\bysame, \emph{Harmonic {A}nalysis for {R}esistance {F}orms}, Japan J. Appl.
  Math. \textbf{204} (2003), 399--444.

\bibitem{RianeDavidM}
N.~Riane and Cl. David, \emph{A spectral study of the {M}inkowski {C}urve,
  hal-01527996}, 2017.

\bibitem{Shima}
T.~Shima, \emph{On eigenvalue problems for the random walks on the
  {S}ierpi\'{n}ski pre-gasket}, Japan J. Indus. Appl. Math. \textbf{8} (1991),
  127--141.

\bibitem{Strichartz1999}
R.~S. Strichartz, \emph{Analysis on fractals}, Notices of the AMS \textbf{46}
  (1999), no.~8, 1199--1208.

\bibitem{StrichartzLivre2006}
\bysame, \emph{Differential equations on fractals, a tutorial}, Princeton
  University Press, 2006.

\bibitem{Strichartz2000}
\bysame, \emph{Taylor approximations on {S}ierpi\'{n}ski {G}asket {T}ype
  {F}ractals}, Journal of Functional Analysis \textbf{174} (2012), 76--127.

\end{thebibliography}
	
\end{document}